\documentclass[12pt,draft,amsmath,amscd,amsbsy,amssymb]{amsart}
\relax

\chardef\bslash=`\\ % p. 424,TeXbook

\makeatletter
\def\verbatim{\interlinepenalty\@M 
\@verbatim
  \leftskip\@totalleftmargin\advance\leftskip2pc
  
\frenchspacing\@vobeyspaces 
\@xverbatim}
\makeatother
\makeatletter
  \def\dgt@k{\dg@DX=-3 
\dg@DY=2 \dg@SIZE=3}
\makeatother

\makeatletter
  \def\dgt@kk{\dg@DX=3 
\dg@DY=-1 
\dg@SIZE=3}%
\makeatother

\theoremstyle{plain}
\newtheorem{thm}{Theorem}[section]
\newtheorem{cor}[thm]{Corollary}
\newtheorem{lem}[thm]{Lemma}
\newtheorem{pro}[thm]{Proposition}

\newtheorem*{A}{Theorem}

\theoremstyle{definition}
\newtheorem{rem}[thm]{Remark}

\numberwithin{equation}{section}

\begin{document}

%%%%%%% Begin Topmatter %%%%%%%%%%

\title[Characterizing the topology 
of
pseudo-boundaries of Euclidean spaces]
{Characterizing the topology 
of
pseudo-boundaries of Euclidean spaces}
\author{Alex 
Chigogidze}
\address{Department of Mathematics and Statistics,
University of 
Saskatche\-wan,
McLean Hall, 106 Wiggins Road, Saskatoon, SK, S7N 
5E6,
Canada}
\email{chigogid@math.usask.ca}
\thanks{The first named author 
was partially supported by
NSERC research 
grant.}

\author{M.~M.~Zarichnyi}
\address{Department of Mechanics and 
Mathematics,
Lviv State University,
Universitetska 1, 290602 Lviv, 
Ukraine}
\email{mzar@litech.lviv.ua}


%%%%%%% End topmatter %%%%%%%%%

\begin{abstract}{We give a 
topological characterization
of the $n$-dimensional pseudo-boundary of 
the
$(2n+1)$-dimensional Euclidean 
space.}
\end{abstract}

\maketitle
\markboth{A.~Chigogidze and 
M.~M.~Zarichnyi}
{Characterizing the topology of
pseudo-boundaries of 
Euclidean spaces}

\section{Introduction}\label{S:intro}
In \cite{gesu} 
Geoghegan and Summerhill constructed
the $n$-dimensional universal 
pseudo-boundary
$\sigma_{n}^{k}$ of the $k$-dimensional Euclidean 
space
${\mathbb R}^{k}$, $0 \leq n \leq k$, $k \geq 1$, as an
${\mathcal 
M}_{n}^{k}$-absorber of ${\mathbb R}^{k}$, where
${\mathcal M}_{n}^{k}$ 
denotes the collection of tame at most
$n$-dimensional compacta in ${\mathbb 
R}^{k}$. In these notes
we consider the space $\sigma_{n}^{2n+1}$.
It has 
been remarked by
several authors that from a certain point of view the 
space
$\sigma_{n}^{2n+1}$ can be considered as the 
$n$-dimensional
counterpart of the
pseudo-boundaries $\sigma$ and $\Sigma$ of 
the Hilbert cube $Q$.
Topological characterizations of the latter spaces have 
been
obtained by Mogilski
\cite{mog}, \cite{bemo}. As for the problem
of 
topological characterization of $\sigma_{n}^{2n+1}$ (see, for instance,
\cite[Problem \# 1017]{west},
\cite[Problem \# 607]{domo}, 
\cite[Conjecture 4.10]{dvmm},
\cite[Question 3]{zar}, \cite[Conjecture 
5.6.9]{book}) we mention
here the following two related results. First of all 
we note that
according to \cite{dvmm}
$\sigma_{n}^{2n+1} \approx 
\sigma_{n}^{k}$ for each
$k \geq 2n+1$. Secondly $\sigma_{n}^{2n+1} \approx 
\Sigma^{n}$
(see \cite[Theorem 7.4]{chikaty}, \cite[Theorem 5.6.10]{book}), 
where
$\Sigma^{n}$ denotes the pseudo-boundary of the 
universal
$n$-dimensional Menger
compactum $\mu^{n}$ \cite{be} constructed in 
\cite{chi1}.

Below (Corollary \ref{C:char}) we give a topological 
characterization of the 
space
$\sigma_{n}^{2n+1}$.

%%%%%%%%%%%%%%%%%%%%%%%%%%%%%%%%%%%%%%%%%%%%%%%%%%%%%%
%%%%%%%%%%%%%%%%%%%%%%%%%%%%%%%%%%%%%%%%%%%%%%%%%%%%

\section{Topological 
characterization of finite-dimensional
absorbing 
sets}\label{S:2}
\subsection{Preliminaries}\label{SS:pre}

All spaces in 
these notes are assumed to be separable
and metrizable.
Maps are assumed to 
be continuous.

Let $n \in \omega$. A subset $A$ of a space $X$ is
said to be 
locally connected in dimension $n$
relative to $X$ (briefly 
$LC^{n}\operatorname{rel.}X$)
if for each $k \leq n+1$, each $x \in X$ and 
each
neighbourhood $U$ of $x$ in $X$ there exists a
neighbourhood $V$ of $x$ 
in $X$ such that every map
$f \colon \partial I^{k} \to V\cap A$ has an 
extension
$F \colon I^{k} \to U \cap A$. A space $X$ is said to be 
locally
connected in dimension $n$ (briefly $LC^{n}$) if $X$ 
is
$LC^{n}\operatorname{rel.}X$. Recall that class of
$LC^{n-1} \cap 
C^{n-1}$-spaces coincides with the class
$AE(n)$ of absolute extensors in 
dimension $n$. Discussion of basic
properties of
$LC^{n}$-spaces can be found 
in \cite{book}.

Recall that for an open cover
${\mathcal U} \in 
\operatorname{cov}(Y)$ of a space $Y$
two maps $f,g \colon X \to Y$ are said 
to be ${\mathcal U}$-close
if for each point $x \in X$ there exists an 
element
$U \in {\mathcal U}$ such that
$f(x), g(x) \in U$.
A space $X$ 
satisfies the discrete $n$-cells
property if the set
\[ \{ f \in 
C(I^{n}\times {\mathbb N}, X) \colon
\{ f(I^{n} \times \{ k\} ) \colon k \in 
{\mathbb N} \}\mbox{ is discrete }\} \]

\noindent is dense in the 
space
$C(I^{n}\times 
{\mathbb 
N}, X)$ equipped
with the limitation topology. The latter
topology on the 
space $C(X,Y)$ of
all continuous maps of $X$ into $Y$ has a 
neighbourhood
base at a point
$f \in C(X,Y)$ consisting of the sets
$\{ g \in 
C(X,Y) \colon g \;\text{is}\;
{\mathcal U}\text{-close to}\; f\}$,
${\mathcal 
U} \in \operatorname{cov}(Y)$.

Proof of our main result is based on the 
following
two statements obtained recently\footnote{The first named 
author
recalls with satisfaction series of very interesting
lectures given by 
S.~Ageev during his stay at the University of
Saskatchewan in May-June of 
1999.} in \cite{ageev}.
\begin{A}[Topological characterization of the 
N\"{o}beling
space \cite{ageev}]\label{T:charactr}
Let $n \geq 0$. Then the 
following conditions are
equivalent for any space 
$X$:
\begin{enumerate}
\item
$X$ is homeomorphic to the $n$-dimensional 
universal
N\"{o}beling space $N_{n}^{2n+1}$.
\item
$X$ is a separable 
completely metrizable space satisfying
the following 
properties:
\begin{itemize}
\item[(a)]
$\dim X = n$.
\item[(b)]
$X \in 
AE(n)$.
\item[(c)]
$X$ has the discrete $n$-cells 
property.
\end{itemize}
\end{enumerate}
\end{A}

\begin{A}[$Z$-set unknotting 
\cite{ageev}]\label{T:unknotting}
Let $n \geq 0$. Then for each open 
cover
${\mathcal U} \in \operatorname{cov}(N_{n}^{2n+1})$ there
exists an 
open cover ${\mathcal V} \in 
\operatorname{cov}(N_{n}^{2n+1})$
such that the 
following property is satisfied:
\begin{itemize}
\item
Every homeomorphism $h 
\colon Z_{1} \to Z_{2}$ between $Z$-sets of
$N_{n}^{2n+1}$ which is 
${\mathcal V}$-close to the inclusion
$Z_{1} \hookrightarrow N_{n}^{2n+1}$ 
can be extended to a
homeomorphism $H \colon N_{n}^{2n+1} \to N_{n}^{2n+1}$ 
which
is ${\mathcal U}$-close to the 
identity
$\operatorname{id}_{N_{n}^{2n+1}}$.
\end{itemize}
\end{A}

\subsection{Uniqueness 
of finite-dimensional absorbing 
sets}\label{SS:uniq}
In this section we 
prove that any two ``absorbing sets" for a class
of finite-dimensional spaces 
are homeomorphic.

Let $\mathcal K$ be a class of spaces that is 
topological,
finitely additive and hereditary with respect to closed 
subspaces.
A space $X$ is {\it strongly} $\mathcal K$-{\it universal}
if, for 
every map $f \colon C
\to X$ from a space $C \in {\mathcal K}$ into $X$, for 
every
closed subspace $D \subseteq C$ such that
$f/D \colon D \to X$ is a 
$Z$-embedding and for every open
cover $\mathcal U \in 
\operatorname{cov}(X)$, there exists
a $Z$-embedding $g \colon C \to X$ such 
that $g/D = f/D$ and
$g$ is $\mathcal U$-close to $f$.

The class consisting 
of countable unions of members of
$\mathcal K$ is denoted by $\mathcal 
K
_{\sigma}$.

Let $n \in \omega$. An $n$-dimensional
separable metrizable 
space $X$ is a $\mathcal
K$-{\it absorbing set}\index{absorbing set} 
if:\\
\begin{itemize}
\item[(a)]
$X \in AE(n)$.
\item[(b)]
$X$ is a countable 
union of strong $Z$-sets.
\item[(c)]
$X \in {\mathcal 
K}_{\sigma}$.
\item[(d)]
$X$ is strongly $\mathcal 
K$-universal.
\end{itemize}

Several examples of $\mathcal K$-absorbing sets 
(for
various classes $\mathcal K$) can be found in
\cite{domo}.

First we 
show that spaces we are interested in can be nicely
embedded into the 
N\"{o}beling space of the same 
dimension.

\begin{pro}\label{P:completion}
Let $n \geq 0$ and $X$ be a 
separable metrizable
$LC^{n-1}\& C^{n-1}$-space satisfying the 
discrete
$n$-cells property. Then $X$ can be embedded into a copy $M$ of 
the
universal $n$-dimensional N\"{o}beling space $N^{2n+1}_{n}$
so that the 
set
$\{ f \in C(I^{n},M) \colon f(I^{n}) \subseteq X\}$
is dense in the space 
$C(I^{n},M)$.
In particular, the following properties are 
satisfied:
\begin{itemize}
\item[(a)]
Every $F_{\sigma}$-subset $F$ of $M$ 
such that
$F \cap X = \emptyset$ is a $Z_{\sigma}$-set 
in
$M$.
\item[(b)]
Every $G_{\delta}$-subspace of $M$, containing $X$, 
is
homeomorphic to $N_{n}^{2n+1}$.
\item[(c)]
If $A$ and $B$ are 
$G_{\delta}$-subsets of $M$ such that
$ X \subseteq A \subseteq B$, then 
$B-A$ is a $Z_{\sigma}$-subset in 
$B$.
\item[(d)]
If $A$ and $B$ are 
$G_{\delta}$-subsets of $M$ such that
$ X \subseteq A \subseteq B$, then the 
inclusion
$A \hookrightarrow B$ is a 
near-homeomorphism.
\end{itemize}
\end{pro}
\begin{proof}
Let 
$\widetilde{X}$
be an $n$-dimensional metrizable compactification of $X$. 
By
\cite[Theorem 2]{banakh}, there exists a $G_{\delta}$-set
$M \subseteq 
\widetilde{X}$,
containing $X$, so that
\begin{itemize}
\item[(1)]
$X$ is 
$\operatorname{LC}^{n-1}\operatorname{rel.}M$;
\item[(2)]
$M \in 
\operatorname{LC}^{n-1}$;
\item[(3)]
For every at most $n$-dimensional Polish 
space $Y$ the set of
all closed embeddings is dense in 
$C(Y,M)$.
\end{itemize}

Let us show that $M \in C^{n-1}$. Indeed, let
$f 
\colon \partial I^{k} \to M$
be a map defined on the boundary $\partial 
I^{k}$ of the
$k$-dimensional disk $I^{k}$, $k \leq n$. According 
to
\cite[Proposition 4.1.7]{book}, there exists an open cover
${\mathcal V} 
\in \operatorname{cov}(M)$ such that the following
condition is 
satisfied:

\begin{itemize}
\item[$(\ast )_{n-1}$]
If a ${\mathcal V}$-close 
to $f$ map $g \colon \partial I^{k} \to M$,
$k \leq n$, has
an extension $G 
\colon I^{k} \to M$, then $f$ also
has an extension $F \colon I^{k} \to 
M$.
\end{itemize}

Since $X$ is 
$\operatorname{LC}^{n-1}\operatorname{rel.}M$, it
follows by \cite[Theorem 
2.8]{tor} that $M-X$ is locally
$n$-negligible in $M$. According to 
\cite[Theorem 2.3]{tor} we
can find a map
$g \colon \partial I^{k} \to X$ 
which is ${\mathcal V}$-close to $f$.
Since $X \in C^{n-1}$, there exists an 
extension
$G \colon I^{k} \to X$ of $g$.
The above stated property $(\ast 
)_{n-1}$ of the cover
${\mathcal V}$ guarantees that
$f$ also has an 
extension $F \colon I^{k} \to M$. This shows that
$M \in C^{n-1}$. Therefore 
$M$ is an $n$-dimensional, separable,
completely metrizable $LC^{n-1}\& 
C^{n-1}$-space satisfying property $(3)$. Topological characterization of 
the N\"{o}beling
space (see Section \ref{S:intro}) implies that $M$ is 
homeomorphic
to $N_{n}^{2n+1}$. The fact that the set the set
$\{ f \in 
C(I^{n}, M) \colon f(I^{n}) \subseteq X\}$
is dense in the space $C(I^{n},M)$ 
follows from \cite[Theorem 
2.8]{book}.

Let $F$ be an $F_{\sigma}$-subset 
of $M$ such that $F \cap X = 
\emptyset$.
Since
\[   \{ f\in C(I^{n},M) 
\colon f(I^{n}) \subseteq X\} \subseteq
\{ f\in C(I^{n},M) \colon f(I^{n}) 
\cap F = \emptyset \} ,\]
it follows that the set
$\{ f\in C(I^{n},M) \colon 
f(I^{n}) \cap F = \emptyset  \}$
is dense in $C(I^{n},M)$.
Consequently,
$F$ 
is a $Z_{\sigma}$-subset of $M$. This proves property (a).

Next observe that 
since
$M$ is homeomorphic to $N_{n}^{2n+1}$ it can be identified with
the 
pseudo-interior $\nu^{n}$ of the universal $n$-dimensional
Menger compactum 
(see \cite[Theorem 5.5.5]{book}). Let $Y$ be a
$G_{\delta}$-subspace of $M$ 
containing $X$. By (a) and
\cite[Proposition 5.7.7]{book}, the inclusion $Y 
\hookrightarrow M$
is a near-homeomorphism. In particular, $Y$ is 
homeomorphic to
$N_{n}^{2n+1}$. This proves (b). Properties (c) and (d) are 
proved 

similarly.
\end{proof}
%%%%%%%%%%%%%%%%%%%%%%%%%%%%%%%%%%%%%%%%%
%%%%%%%%%%%%%%%%%%%%%%%%%%%%%%%%%%%%%%%%

\begin{rem}
An 
$n$-dimensional $\mathcal{K}$-absorbing set is called {\it 
representable\/} 
in $\mathbb{R}^k$ \cite{domo} if there exists an 
embedding
$i\colon 
M\to\mathbb{R}^k$ such that the set $\mathbb{R}^k-M$ is 
locally 

$n$-negigible in $\mathbb{R}^k$. 
\begin{itemize}
\item[ ]
{\em Every 
$n$-dimensional $\mathcal{K}$-absorbing set is 
representable in 
$\mathbb{R}^{2n+1}$.}
\end{itemize}
\begin{proof} It is showm in the proof of 
Proposition 
\ref{P:completion} that 
there exists an embedding of 
$\mathcal{K}$-absorbing set $M$ into 
$N^{2n+1}_{n}$ 
with locally 
$n$-negligible complement of the image. Now observe that
the complement 
${\mathbb R}^{2n+1} - N_{n}^{2n+1}$ as
a $\sigma Z_{n}$-set \cite{book} in 
${\mathbb R}^{2n+1}$ is locally
$n$-negligible in ${\mathbb R}^{2n+1}$. This 
obviously implies
that the complement
${\mathbb R}^{2n+1}$ is also 
$n$-negligible in
${\mathbb R}^{2n+1}$ as required.
\end{proof}
The above 
statement provides an affirmative solution of Problem
$555$ from 
\cite{domo}.
\end{rem}
%%%%%%%%%%%%%%%%%%%%%%%%%%%%%%%%%%%%%%%

\begin{lem}\label{L:compact}
Let 
$X$ be an at most $n$-dimensional separable metrizable
$LC^{n-1}$-space. 
If
$X = \cup\{ X_{i} \colon i \in \omega\}$, where each
$X_{i}$ is a strong 
$Z$-set in $X$, then each compact
subset of $X$ is a strong $Z$-set in 
$X$.
\end{lem}
\begin{proof}
Let $\widetilde{X}$ be an $n$-dimensional 
separable completely
metrizable space
containing $X$ as a subspace in such a 
way that $X$ 
is
$\operatorname{LC}^{n-1}\operatorname{rel.}\widetilde{X}$
(see 
\cite[Proposition 2.8]{domar}.
As in the proof
of Proposition 
\ref{P:completion}, we conclude that
\begin{itemize}
\item[$(\ast )$]
the 
set
$\{ f \in C(I^{n},\widetilde{X}) \colon f(I^{n}) \subseteq X \}$
is dense 
in $C(I^{n},\widetilde{X})$.
\end{itemize}
Next we need the following 
observation.

{\bf Claim}. {\em A compact subset $K$ of $X$ is a $Z$-set in 
$X$
if and only if $K$ is a $Z$-set in $\widetilde{X}$.}

{\em Proof of 
Claim}. First let $K$ be a $Z$-set in $X$. Consider a 
map
$f \colon I^{n} 
\to \widetilde{X}$ and open covers
${\mathcal U}, {\mathcal V} \in 
\operatorname{cov}(\widetilde{X})$
such that $\operatorname{St}({\mathcal 
V})$ refines ${\mathcal U}$.
By $(\ast )$, there exists a
${\mathcal 
V}$-close to $f$ map $g \in C(I^{n},\widetilde{X})$ such
that $g(I^{n}) 
\subseteq X$. Since $K$ is a $Z$-set in $X$,
there exists
a ${\mathcal 
V}$-close to $g$ map $h \colon I^{n} \to X$ such that
$h(I^{n}) \cap K = 
\emptyset$. Since $h$ is ${\mathcal U}$-close to
$f$, it follows that $K$ is 
a $Z$-set in $\widetilde{X}$.

Conversely, let $K$ be a $Z$-set in 
$\widetilde{X}$. Consider a map
$f \colon I^{n} \to X$ and open 
covers
${\mathcal U}, {\mathcal V} \in \operatorname{cov}(X)$ so 
that
$\operatorname{St}({\mathcal V})$ refines ${\mathcal U}$.
For each $V 
\in {\mathcal V}$ choose an open subset
$\widetilde{V} \subseteq 
\widetilde{X}$
such that $V = \widetilde{V} \cap X$. It is easy to see 
that
$K$ is a $Z$-set in
$Y = \cup\{ \widetilde{V} \colon V \in {\mathcal 
V}\}$.
Consequently there exists
a $\widetilde{{\mathcal V}}$-close to $f$ 
map
$g \colon I^{n} \to Y$ such that
$g(I^{n}) \cap K = \emptyset$, 
where
$\widetilde{{\mathcal V}} = \{ \widetilde{V} \colon V
\in {\mathcal 
V}\} \in \operatorname{cov}(Y)$. Let $G$ be an
open subsets of $Y$ such that 
$K \cap G = \emptyset$
and $g(I^{n}) \subseteq G$. By $(\ast )$, there exists 
a map
$h \colon I^{n} \to X$
which is
$\widetilde{{\mathcal V}} \wedge \{ G, 
Y - g(I^{n})\}$-close
to $g$. Obviously, $h$ is ${\mathcal U}$-close to $f$ 
and
$h(I^{n}) \cap K = \emptyset$. This shows that $K$ is a
$Z$-set in $X$ 
and completes the proof of claim.

We continue the proof of Lemma 
\ref{L:compact}. Let $K$ be
a compact subset of $X$. Clearly $K \cap X_{i}$ 
is a compact
$Z$-set in $X$ for each $i \in \omega$. By the above Claim,
$K 
\cap X_{i}$ is a $Z$-set in $\widetilde{X}$. This means that
the set
$\{ f 
\in C(I^{n}, \widetilde{X}) \colon
f(I^{n}) \cap (X_{i}\cap K ) = 
\emptyset\}$ is open and dense
in the space $C(I^{n},\widetilde{X})$. Since 
$\widetilde{X}$
is completely metrizable, the space $C(I^{n}, \widetilde{X})$ 
has
the Baire property (see, for instance, \cite[Proposition 
2.1.7]{book}) 
and
consequently, the set
\begin{multline*}
\{ f \in C(I^{n},\widetilde{X}) 
\colon f(I^{n}) \cap K = \emptyset \} 
= \\
\{ f \in C(I^{n},\widetilde{X}) 
\colon f(I^{n}) \cap
\left( \cup\{ X_{i}\cap K \colon i \in \omega\} \right) 
= \emptyset 
\} =\\
\bigcap\left\{ \{ f \in C(I^{n},\widetilde{X}) \colon 
f(I^{n})
\cap (X_{i} \cap K) = \emptyset\} \colon i \in \omega 
\right\}
\end{multline*}
is also dense in the space $C(I^{n},\widetilde{X})$. 
This simply 
means
that $K$ is a $Z$-set in $\widetilde{X}$. By the above 
Claim,
we conclude that $K$ is a $Z$-set in $X$ as well.
\end{proof}

Proof 
of the following statement uses Lemma \ref{L:compact} and
follows verbatim 
the proof of
\cite[Lemma 1.9]{domo}.

\begin{pro}\label{P:ok}
Let $X$ be 
an at most $n$-dimensional separable metrizable
$LC^{n-1}$-space. If
$X = 
\cup\{ X_{i} \colon i \in \omega\}$, where each
$X_{i}$ is a strong $Z$-set 
in $X$, then
$X$ satisfies the discrete $n$-cells property.
\end{pro}

Now 
we are in position to prove the uniqueness 
result.

\begin{thm}\label{T:main}
Let $n \geq 0$ and $\mathcal K$ be a class 
of spaces that is
topological, finitely additive and hereditary with respect 
to
closed subspaces. Then any two ${\mathcal K}(n)$-absorbing sets
are 
homeomorphic.
\end{thm}
\begin{proof}
Let $X$ and $Y$ be ${\mathcal 
K}(n)$-absorbing sets.
Proposition
\ref{P:ok} guarantees that $X, Y \in 
n$-$\operatorname{SDAP}$.
Embed $X$ and $Y$ into a copy $M$ of the 
universal
$n$-dimensional
N\"{o}beling space $N_{n}^{2n+1}$ in such a way 
that
properties (a)--(d) of Proposition \ref{P:completion}
are satisfied. The 
rest of the proof follows the argument
presented in the proof 
of
\cite[Theorem 3.1]{bemo} (use the $Z$-set Unknotting
Theorem for 
$N_{n}^{2n+1}$ instead of the $Z$-set
unknotting theorem for $\ell_{2}$ at 
the appropriate 
place).
\end{proof}

%%%%%%%%%%%%%%%%%%%%%%%%%%%%%%%%%%%%%%%%%%%%%%%%%%%%%
%%%%%%%%%%%%%%%%%%%%%%%%%%%%%%%%%%%%%%%%%%%%%%%%%%%%%
\subsection{Characterization 
of 
$\sigma_{n}^{2n+1}$}\label{SS:charact}
In order to obtain a topological 
characterization of a
${\mathcal K}(n)$-absorbing set Theorem \ref{T:main} 
must be
combined with the corresponding existence result. In other words, 
we
need to know that there exists a ${\mathcal K}(n)$-absorbing set.
For 
certain choices of ${\mathcal K}$ it is even possible to
explicitly construct 
corresponding absorbing sets.

Let us recall that for each space
$X$ and for 
each ordinal
$\alpha < {\omega}_1$, we can define
two classes of subspaces of 
$X$ -- the {\it additive Borelian
class} $\alpha$, $\mathcal A 
_{\alpha}(X)$,
and the {\it multiplicative Borelian class} 
$\alpha$,
$\mathcal M _{\alpha}(X)$, -- as follows: $\mathcal A _0
(X)$ is 
the collection of all open subsets of $X$ and
$\mathcal M _0 (X)$ is the 
collection of all closed
subsets of $X$. Assuming that for each 
ordinal
$\beta < \alpha$, where $\alpha < {\omega}_1$,
the classes $\mathcal 
A _{\beta}$ and $\mathcal M _{\beta}$
have already been constructed, we 
proceed as
follows: the class $\mathcal A _{\alpha}$ consists of
countable 
unions of elements of $\cup \{ \mathcal M
_{\beta} \colon \beta < \alpha \}$ 
and the class
$\mathcal M _{\alpha}$ consists of countable
intersections of 
elements of
$\cup \{ \mathcal A _{\beta} \colon \beta < \alpha \}$.

Further, 
let $\alpha < {\omega}_1$ and $X$ be a  separable
metrizable space. We say 
that $X$
belongs to the {\it absolute additive Borelian class}
$\mathcal A 
_{\alpha}$, if for any embedding $i \colon
X \to Y$ into any separable 
metrizable space $Y$, we have
$i(X) \in \mathcal A _{\alpha}(Y)$. 
Similarly,
$X$ belongs to the {\it absolute multiplicative Borelian
class} 
$\mathcal M _{\alpha}$ if for any
embedding $i \colon X \to Y$ into any 
separable metrizable
space $Y$, we have $i(X) \in \mathcal M
_{\alpha}(Y)$. 
It is well-known that:
(a) $X \in \mathcal A _{\alpha}$, $\alpha \geq 2$, if 
and
only if $X \in \mathcal A _{\alpha}(l_2 )$ and
(b) $X \in \mathcal M 
_{\alpha}$, $\alpha \geq 1$, if and
only if $X \in \mathcal M _{\alpha}(l_2 
)$.

Obviously, $\mathcal A _0 = \emptyset$ and $\mathcal M _0$
coincides 
with the class of all metrizable
compacta. Further,
$\mathcal A _1 = \{ 
\sigma$-compact spaces$\}$,
$\mathcal M _1 = \{ \mbox{Polish
spaces} \}$, 
etc.

The existence problem for these classes of spaces is solved in
the 
following statement \cite[Theorem 5.7.21]{book},
\cite[Theorem 
2.5]{zar}.

\begin{thm}\label{T:ex}
Let $n \in \omega$ and $1 \leq \alpha < 
\omega_{1}$. Then there exist
an ${\mathcal A}_{\alpha}(n)$-absorbing set 
$\Lambda_{\alpha}(n)$
and ${\mathcal M}_{\alpha}(n)$-absorbing set 
$\Omega_{\alpha}(n)$.
\end{thm}

Theorems \ref{T:main} and \ref{T:ex} imply 
the following
characterization result.

\begin{thm}\label{T:gen}
Let $X$ be 
an $n$-dimensional, $n \geq 0$, separable metrizable
$AE(n)$-space and $1 
\leq \alpha <\omega_{1}$. Then $X$ is
homeomorphic to $\Omega_{\alpha}(n)$ 
(respectively,
$\Lambda_{\alpha}(n)$) if and only if the following two 
conditions 
are
satisfied:
\begin{itemize}
\item[(i)]
$X = \cup\{ X_{i} 
\colon i \in \omega\}$, where each
$X_{i} \in {\mathcal M}_{\alpha}$ 
(respectively,
$X_{i} \in {\mathcal A}_{i}$) and $X_{i}$ is a strong $Z$-set 
in $X$,
\item[(ii)]
$X$ is strongly ${\mathcal M}_{\alpha}(n)$-universal 
(respectively,
${\mathcal 
A}_{\alpha}(n)$-universal.
\end{itemize}
\end{thm}

In particular ($\alpha = 
1$), we obtain a topological
characterization of 
$\sigma_{n}^{2n+1}$.

\begin{cor}\label{C:char}
Let $X$ be an 
$n$-dimensional, $n \geq 0$, $\sigma$-compact 
metrizable
$AE(n)$-space. Then 
$X$ is homeomorphic to $\sigma_{n}^{2n+1}$
if and only if the following 
conditions are satisfied:
\begin{itemize}
\item[(i)]
$X$ has the discrete 
$n$-cells property,
\item[(ii)]
$X$ is strongly ${\mathcal 
A}_{1}(n)$-universal.
\end{itemize}
\end{cor}
%%%%%%%%%%%%%%%%%%%%%%%%%%%%%%%%%%%%%%%%%%%%%%%%%%%%%%%%

\end{document}